# Bertrand Partner $D$-Curves in Euclidean 3-space $E^3$

**Mustafa Kazaz[a], H. Hüseyin Uğurlu[b], Mehmet Önder[a], Seda Oral[a]**

[a]*Celal Bayar University, Department of Mathematics, Faculty of Arts and Sciences, , Manisa, Turkey.*
E-mails:mustafa.kazaz@bayar.edu.tr, mehmet.onder@bayar.edu.tr
[b]*Gazi University, Gazi Faculty of Education, Department of Secondary Education Science and Mathematics Teaching, Mathematics Teaching Program, Ankara, Turkey.* E-mail: hugurlu@gazi.edu.tr

**Abstract**
In this paper we consider the idea of Bertrand curves for curves lying on surfaces and by considering the Darboux frames of them we define these curves as Bertrand $D$-curves and give the characterizations for these curves. We also find the relations between the geodesic curvatures, the normal curvatures and the geodesic torsions of these associated curves. Furthermore, we show that the definition and the characterizations of Bertrand $D$-curves include those of Bertrand curves in some special cases.

## 1. Introduction

Bertrand curves are one of the associated curve pairs for which at the corresponding points of the curves one of the Frenet vectors of a curve coincides with the one of the Frenet vectors of the other curve. These special curves are very interesting and an important problem of the fundamental theory and the characterizations of space curve and are characterized as a kind of corresponding relation between the two curves such that the curves have the common principal normal i.e., the Bertrand curve is a curve which shares the normal line with another curve. These curves have an important role in the theory of curves and surfaces. Hereby, from the past to today, a lot of mathematicians have studied on Bertrand curves in different areas[1,2,3,5,13,14]. In [8], Izumiya and Takeuchi have studied cylindrical helices and Bertrand curves from the view point as curves on ruled surfaces. They have shown that cylindrical helices can be constructed from plane curves and Bertrand curves can be constructed from spherical curves. Also, they have studied generic properties of cylindrical helices and Bertrand curves as applications of singularity theory for plane curves and spherical curves[9]. In [4], Gluck has investigated the Bertrand curves in $n$-dimensional Euclidean space $E^n$. The corresponding characterizations of the Bertrand curves in $n$-dimensional Lorentzian space have been given by Ekmekçi and Ilarslan[3]. They have shown that the distance between corresponding points of the Bertrand couple curves and the angle between the tangent vector fields of these points is constant.

Furthermore, by considering the Frenet frame of the ruled surfaces, Ravani and Ku extended the notion of Bertrand curve to the ruled surfaces and named as Bertrand offsets[12]. The corresponding characterizations of the Bertrand offsets of timelike ruled surface were given by Kurnaz[10].

In this paper, we consider the notion of the Bertrand curve for the curves lying on the surfaces. We call these new associated curves as Bertrand $D$-curves and by using the Darboux frame of the curves we give the definition and the characterizations of these curves.

## 2. Darboux Frame of a Curve Lying on a Surface

Let $S = S(u,v)$ be an oriented surface in three-dimensional Euclidean space $E^3$ and let consider a curve $x(s)$ lying on $S$ fully. Since the curve $x(s)$ is also in space, there exists Frenet frame $\{T, N, B\}$ at each points of the curve where $T$ is unit tangent vector, $N$ is principal normal vector and $B$ is binormal vector, respectively. The Frenet equations of the curve $x(s)$ is given by

$$T' = \kappa N$$
$$N' = -\kappa T + \tau B$$
$$B' = -\tau N$$

where $\kappa$ and $\tau$ are curvature and torsion of the curve $x(s)$, respectively.

Since the curve $x(s)$ lies on the surface $S$ there exists another frame of the curve $x(s)$ which is called Darboux frame and denoted by $\{T, g, n\}$. In this frame $T$ is the unit tangent of the curve, $n$ is the unit normal of the surface $S$ and $g$ is a unit vector given by $g = n \times T$. Since the unit tangent $T$ is common in both Frenet frame and Darboux frame, the vectors $N$, $B$, $g$ and $n$ lie on the same plane. So that the relations between these frames can be given as follows

$$\begin{bmatrix} T \\ g \\ n \end{bmatrix} = \begin{bmatrix} 1 & 0 & 0 \\ 0 & \cos\varphi & \sin\varphi \\ 0 & -\sin\varphi & \cos\varphi \end{bmatrix} \begin{bmatrix} T \\ N \\ B \end{bmatrix}$$

where $\varphi$ is the angle between the vectors $g$ and $N$. The derivative formulae of the Darboux frame is

$$\begin{bmatrix} \dot{T} \\ \dot{g} \\ \dot{n} \end{bmatrix} = \begin{bmatrix} 0 & k_g & k_n \\ -k_g & 0 & \tau_g \\ -k_n & -\tau_g & 0 \end{bmatrix} \begin{bmatrix} T \\ g \\ n \end{bmatrix} \quad (1)$$

where $k_g$, $k_n$ and $\tau_g$ are called the geodesic curvature, the normal curvature and the geodesic torsion, respectively. Here and in the following, we use "dot" to denote the derivative with respect to the arc length parameter of a curve.

The relations between geodesic curvature, normal curvature, geodesic torsion and $\kappa$, $\tau$ are given as follows

$$k_g = \kappa \cos\varphi, \quad k_n = \kappa \sin\varphi, \quad \tau_g = \tau + \frac{d\varphi}{ds}. \quad (2)$$

Furthermore, the geodesic curvature $k_g$ and geodesic torsion $\tau_g$ of the curve $x(s)$ can be calculated as follows

$$k_g = \left\langle \frac{dx}{ds}, \frac{d^2x}{ds^2} \times n \right\rangle, \quad \tau_g = \left\langle \frac{dx}{ds}, n \times \frac{dn}{ds} \right\rangle \quad (3)$$

In the differential geometry of surfaces, for a curve $x(s)$ lying on a surface $S$ the followings are well-known

**i)** $x(s)$ is a geodesic curve $\Leftrightarrow k_g = 0$,

**ii)** $x(s)$ is an asymptotic line $\Leftrightarrow k_n = 0$,

**iii)** $x(s)$ is a principal line $\Leftrightarrow \tau_g = 0$ [11].

Furthermore, if $(k_g)_1$ and $(k_g)_2$ are geodesic curvatures of the parametric curves of $S(u,v)$ and $\sigma$ is the angle between the curve $x(s)$ and the parametric curve $v = constant$ of $S(u,v)$ then we have

$$k_g = \frac{d\sigma}{ds} + (k_g)_1 \cos\sigma + (k_g)_2 \sin\sigma$$

which is known Liouville's formula (See [13]).

## 3. Bertrand $D$-Curves in Euclidean 3-space $E^3$

In this section, by considering the Darboux frame, we define Bertrand $D$-curves and give the characterizations of these curves.

**Definition 1.** Let $S$ and $S_1$ be oriented surfaces in three-dimensional Euclidean space $E^3$ and let consider the arc-length parameter curves $x(s)$ and $x_1(s_1)$ lying fully on $S$ and $S_1$, respectively. Denote the Darboux frames of $x(s)$ and $x_1(s_1)$ by $\{T, g, n\}$ and $\{T_1, g_1, n_1\}$, respectively. If there exists a corresponding relationship between the curves $x$ and $x_1$ such that, at the corresponding points of the curves, the Darboux frame element $g$ of $x$ coincides with the Darboux frame element $g_1$ of $x_1$, then $x$ is called a Bertrand $D$-curve, and $x_1$ is a Bertrand partner $D$-curve of $x$. Then, the pair $\{x, x_1\}$ is said to be a Bertrand $D$-pair. If there exist such curves lying on the oriented surfaces $S$ and $S_1$, respectively, we call the pair $\{S, S_1\}$ as Bertrand pair surfaces.

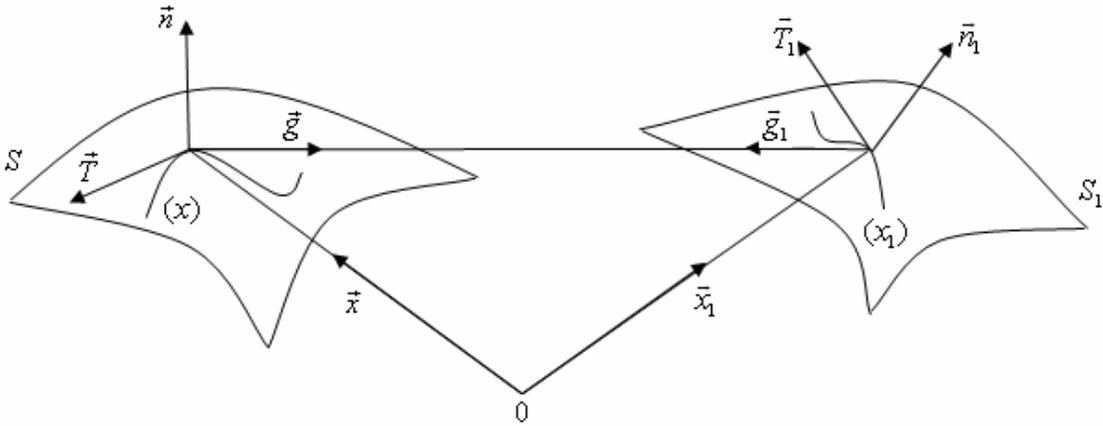

**Fig. 1** Bertrand partner $D$-curves

**Theorem 1.** Let $S$ be an oriented surface and $x(s)$ be a Bertrand $D$-curve in $E^3$ with arc length parameter $s$ fully lying on $S$. If $S_1$ is another oriented surface and $x_1(s_1)$ is a curve with arc length parameter $s_1$ fully lying on $S_1$, then $x_1(s_1)$ is Bertrand partner $D$-curve of $x(s)$ if and only if the normal curvature $k_n$ of $x(s)$ and the geodesic curvature $k_{g_1}$, the normal curvature $k_{n_1}$ and the geodesic torsion $\tau_{g_1}$ of $x_1(s_1)$ satisfy the following equation

$$-\lambda \dot{\tau}_{g_1} = \left( \frac{(1-\lambda k_{g_1})^2 + \lambda^2 \tau_{g_1}^2}{(1-\lambda k_{g_1})} \right)\left( k_{n_1} - k_n \frac{1-\lambda k_{g_1}}{\cos\theta} \right) + \frac{\lambda^2 \tau_{g_1} \dot{k}_{g_1}}{1-\lambda k_{g_1}}$$

for some nonzero constants $\lambda$, where $\theta$ is the angle between the tangent vectors $T$ and $T_1$ at the corresponding points of $x$ and $x_1$.

**Proof:** Suppose that $S$ is an oriented surface and $x(s)$ is a Bertrand $D$-curve fully lying on $S$. Denote the Darboux frames of $x(s)$ and $x_1(s_1)$ by $\{T, g, n\}$ and $\{T_1, g_1, n_1\}$, respectively. Then by the definition we can assume that

$$x(s_1) = x_1(s_1) + \lambda(s_1) g_1(s_1), \qquad (4)$$

for some function $\lambda(s_1)$. By taking derivative of (4) with respect to $s_1$ and applying the Darboux formulas (1) we have

$$T\frac{ds}{ds_1} = (1-\lambda k_{g_1})T_1 + \dot{\lambda}g_1 + \lambda\tau_{g_1}n_1 \tag{5}$$

Since the direction of $g_1$ coincides with the direction of $g$, i.e., the tangent vector $T$ of the curve lies on the plane spanned by the vectors $T_1$ and $n_1$, we get

$$\dot{\lambda}(s_1) = 0.$$

This means that $\lambda$ is a nonzero constant. Thus, the equality (5) can be written as follows

$$T\frac{ds}{ds_1} = (1-\lambda k_{g_1})T_1 + \lambda\tau_{g_1}n_1. \tag{6}$$

Furthermore, we have

$$T = \cos\theta T_1 - \sin\theta n_1, \tag{7}$$

where $\theta$ is the angle between the tangent vectors $T$ and $T_1$ at the corresponding points of $x$ and $x_1$. By differentiating this last equation with respect to $s_1$, we get

$$(k_g g + k_n n)\frac{ds}{ds_1} = -(\dot{\theta}+k_{n_1})\sin\theta T_1 + (k_{g_1}\cos\theta + \tau_{g_1}\sin\theta)g_1 + (-\dot{\theta}+k_{n_1})\cos\theta n_1. \tag{8}$$

From this equation and the fact that

$$n = \sin\theta T_1 + \cos\theta n_1, \tag{9}$$

we get

$$(k_n\sin\theta T_1 + k_g g + k_n\cos\theta n_1)\frac{ds}{ds_1} = (-\dot{\theta}+k_{n_1})\sin\theta T_1 + (k_{g_1}\cos\theta + \tau_{g_1}\sin\theta)g_1 \tag{10}$$
$$+(-\dot{\theta}+k_{n_1})\cos\theta n_1$$

Since the direction of $g_1$ is coincident with $g$ we have

$$\dot{\theta} = k_{n_1} - k_n\frac{ds}{ds_1}. \tag{11}$$

From (6) and (7) and notice that $T_1$ is orthogonal to $g_1$ we obtain

$$\frac{ds}{ds_1} = \frac{1-\lambda k_{g_1}}{\cos\theta} = \frac{-\lambda\tau_{g_1}}{\sin\theta}. \tag{12}$$

Equality (12) gives us

$$\tan\theta = \frac{-\lambda\tau_{g_1}}{1-\lambda k_{g_1}}. \tag{13}$$

By taking the derivative of this equation and applying (11) we get

$$-\lambda\dot{\tau}_{g_1} = \left(\frac{(1-\lambda k_{g_1})^2 + \lambda^2\tau_{g_1}^2}{(1-\lambda k_{g_1})}\right)\left(k_{n_1} - k_n\frac{1-\lambda k_{g_1}}{\cos\theta}\right) + \frac{\lambda^2\tau_{g_1}\dot{k}_{g_1}}{1-\lambda k_{g_1}}, \tag{14}$$

that is desired.

Conversely, assume that the equation (14) holds for some nonzero constants $\lambda$. Then by using (12) and (13), (14) gives us

$$k_n\left(\frac{ds}{ds_1}\right)^3 = \lambda\dot{\tau}_{g_1}(1-\lambda k_{g_1}) + \lambda^2\tau_{g_1}\dot{k}_{g_1} + \left((1-\lambda k_{g_1})^2 + \lambda^2\tau_{g_1}^2\right)k_{n_1} \tag{15}$$

Let define a curve

$$x(s_1) = x_1(s_1) + \lambda(s_1)g_1(s_1). \tag{16}$$

We will prove that $x$ is a Bertrand $D$-curve and $x_1$ is the Bertrand partner $D$-curve of $x$. By taking the derivative of (16) with respect to $s_1$ twice, we get

$$T\frac{ds}{ds_1} = (1-\lambda k_{g_1})T_1 + \lambda \tau_{g_1} n_1, \qquad (17)$$

and

$$(k_g g + k_n n)\left(\frac{ds}{ds_1}\right)^2 + T\frac{d^2 s}{ds_1^2} = -\lambda(\dot{k}_{g_1} + \tau_{g_1} k_{n_1})T_1 + \left((1-\lambda k_{g_1})k_{g_1} - \lambda \tau_{g_1}^2\right)g_1 \\ + \left((1-\lambda k_{g_1})k_{n_1} + \lambda \dot{\tau}_{g_1}\right)n_1 \qquad (18)$$

respectively. Taking the cross product of (17) with (18) we have

$$\left[k_g n - k_n g\right]\left(\frac{ds}{ds_1}\right)^3 = \left[-\lambda \tau_{g_1} k_{g_1}(1-\lambda k_{g_1}) + \lambda^2 \tau_{g_1}^3\right]T_1 \\ - \left[(1-\lambda k_{g_1})^2 k_{n_1} + \lambda \dot{\tau}_{g_1}(1-\lambda k_{g_1}) + \lambda^2 \tau_{g_1} \dot{k}_{g_1} + \lambda^2 \tau_{g_1}^2 k_{n_1}\right]g_1 \\ + \left[k_{g_1}(1-\lambda k_{g_1})^2 - \lambda \tau_{g_1}^2(1-\lambda k_{g_1})\right]n_1 \qquad (19)$$

By substituting (15) in (19) we get

$$\left[k_g n - k_n g\right]\left(\frac{ds}{ds_1}\right)^3 = \left(-\lambda \tau_{g_1} k_{g_1}(1-\lambda k_{g_1}) + \lambda^2 \tau_{g_1}^3\right)T_1 - k_g\left(\frac{ds}{ds_1}\right)^3 g_1 \\ + \left(k_{g_1}(1-\lambda k_{g_1})^2 - \lambda \tau_{g_1}^2(1-\lambda k_{g_1})\right)n_1 \qquad (20)$$

Taking the cross product of (17) with (20) we have

$$\left[-k_g g - k_n n\right]\left(\frac{ds}{ds_1}\right)^4 = \lambda k_n \tau_{g_1}\left(\frac{ds}{ds_1}\right)^3 T_1 + \left((1-\lambda k_{g_1})^2 + \lambda^2 \tau_{g_1}^2\right)\left(\lambda \tau_{g_1}^2 - k_{g_1}(1-\lambda k_{g_1})\right)g_1 \\ - k_n(1-\lambda k_{g_1})\left(\frac{ds}{ds_1}\right)^3 n_1 \qquad (21)$$

From (20) and (21) we have

$$-(k_g^2 + k_n^2)\left(\frac{ds}{ds_1}\right)^4 n = \left[\lambda k_g k_{g_1}\tau_{g_1}(1-\lambda k_{g_1})\frac{ds}{ds_1} - \lambda^2 k_g \tau_{g_1}^3 \frac{ds}{ds_1} + \lambda \tau_{g_1} k_n^2\left(\frac{ds}{ds_1}\right)^3\right]T_1 \\ + k_n\left(\frac{ds}{ds_1}\right)^2\left[k_g\left(\frac{ds}{ds_1}\right)^2 + \lambda \tau_{g_1}^2 - k_{g_1}(1-\lambda k_{g_1})\right]g_1 \\ + \left[-k_g k_{g_1}(1-\lambda k_{g_1})^2 \frac{ds}{ds_1} + \lambda \tau_{g_1}^2 k_g(1-\lambda k_{g_1})\frac{ds}{ds_1} - k_n^2(1-\lambda k_{g_1})\left(\frac{ds}{ds_1}\right)^3\right]n_1 \qquad (22)$$

Furthermore, from (17) and (20) we get

$$\begin{cases} \left(\frac{ds}{ds_1}\right)^2 = (1-\lambda k_{g_1})^2 + \lambda^2 \tau_{g_1}^2, \\ k_g\left(\frac{ds}{ds_1}\right)^2 = k_{g_1}(1-\lambda k_{g_1}) - \lambda \tau_{g_1}^2 \end{cases} \qquad (23)$$

respectively. Substituting (23) in (22) we obtain

$$-(k_g^2 + k_n^2)\left(\frac{ds}{ds_1}\right)^4 n = \left[\lambda k_g k_{g_1}\tau_{g_1}(1-\lambda k_{g_1})\frac{ds}{ds_1} - \lambda^2 k_g \tau_{g_1}^3 \frac{ds}{ds_1} + \lambda \tau_{g_1} k_n^2\left(\frac{ds}{ds_1}\right)^3\right]T_1 \\ + \left[-k_g^2(1-\lambda k_{g_1})^2 \frac{ds}{ds_1} + \lambda \tau_{g_1}^2 k_g(1-\lambda k_{g_1})\frac{ds}{ds_1} - k_n^2(1-\lambda k_{g_1})\left(\frac{ds}{ds_1}\right)^3\right]n_1 \qquad (24)$$

Equality (17) and (24) shows that the vectors $\vec{T}$ and $\vec{n}$ lie on the plane $sp\{\vec{T}_1,\vec{n}_1\}$. So, at the corresponding points of the curves, the Darboux frame element $\vec{g}$ of $x$ coincides with the Darboux frame element $\vec{g}_1$ of $x_1$, i.e, the curves $x$ and $x_1$ are Bertrand $D$-pair curves.

Let now give the characterizations of Bertrand partner $D$-curves in some special cases. Assume that $x(s)$ be an asymptotic Bertrand $D$-curve. Then, from (14) we have the following special cases:

**i)** Consider that $x_1(s_1)$ is a geodesic curve. Then $x_1(s_1)$ is Bertrand partner $D$-curve of $x(s)$ if and only if the following equation holds,
$$\lambda \dot{\tau}_{g_1} = -k_{n_1}(1+\lambda^2 \tau_{g_1}^2)$$

**ii)** Assume that $x_1(s_1)$ is also an asymptotic line. Then $x_1(s_1)$ is Bertrand partner $D$-curve of $x(s)$ if and only if the geodesic torsion $\tau_{g_1}$ of $x_1(s_1)$ satisfies the following equation,
$$\dot{\tau}_{g_1} = -\frac{\lambda \tau_{g_1} \dot{k}_{g_1}}{1-\lambda k_{g_1}}.$$

In this case, the Frenet frame of the curve $x_1(s_1)$ coincides with its Darboux frame. From (2) we have $k_{g_1} = \kappa_1$ and $\tau_{g_1} = \tau_1$. So, the Bertrand partner $D$-curves become the Bertrand partner curves, i.e., if both $x(s)$ and $x_1(s_1)$ are asymptotic lines then, the definition and the characterizations of the Bertrand partner $D$-curves involve those of the Bertrand partner curves in Euclidean 3-space. Then, a new characterization of Bertrand curves can be given as follows

**Corollary 1.** *Let $x(s)$ be a Bertrand curve with arclength parameter $s$. Then the curve $x_1(s_1)$ is a Bertrand partner curve of $x(s)$ if and only if the curvature $\kappa_1$ and the torsion $\tau_1$ of $x_1(s_1)$ satisfy the following equation*
$$\dot{\tau}_1 = -\frac{\lambda \dot{\kappa}_1 \tau_1}{1-\lambda \kappa_1}$$
*for some nonzero constants $\lambda$.*

**iii)** If $x_1(s_1)$ is a principal line then $x_1(s_1)$ is Bertrand partner $D$-curve of $x(s)$ if and only if the geodesic curvature $k_{g_1}$ and the geodesic torsion $\tau_{g_1}$ of $x_1(s_1)$ satisfy the following equality,
$$k_{n_1}(1-\lambda k_{g_1}) = 0.$$

**Theorem 2.** *Let the pair $\{x, x_1\}$ be a Bertrand $D$-pair. Then the relation between geodesic curvature $k_g$, geodesic torsion $\tau_g$ of $x(s)$ and the geodesic curvature $k_{g_1}$, the geodesic torsion $\tau_{g_1}$ of $x_1(s_1)$ is given as follows*
$$k_g - k_{g_1} = \lambda(k_g k_{g_1} - \tau_g \tau_{g_1}).$$

**Proof:** Let $x(s)$ be a Bertrand $D$-curve and $x_1(s_1)$ be a Bertrand partner $D$-curve of $x(s)$. Then from (16) we can write
$$x_1(s_1) = x(s_1) - \lambda g_1(s_1) \tag{23}$$
for some constants $\lambda$. By differentiating (23) with respect to $s_1$ we have
$$\vec{T}_1 = (1+\lambda k_g)\frac{ds}{ds_1}\vec{T} - \lambda \tau_g \frac{ds}{ds_1}\vec{n} \tag{24}$$

By the definition we have
$$\vec{T}_1 = \cos\theta \vec{T} + \sin\theta \vec{n} \tag{25}$$
From (24) and (25) we obtain
$$\cos\theta = (1\lambda k_{g_1})\frac{ds_1}{ds}, \quad \sin\theta = -\lambda \tau_{g_1}\frac{ds_1}{ds} \tag{26}$$
Using (12) and (26) it is easily seen that
$$k_g - k_{g_1} = \lambda(k_g k_{g_1} - \tau_g \tau_{g_1}).$$

From Theorem 2, we obtain the following special cases.
Let the pair $\{x, x_1\}$ be a Bertrand $D$-pair. Then,

i) if one of the curves $x$ and $x_1$ is a principal line, then the relation between the geodesic curvatures $k_g$ and $k_{g_1}$ is
$$k_g - k_{g_1} = \lambda k_g k_{g_1}$$

ii) if $x_1$ is a geodesic curve, then the geodesic curvature of the curve $x$ is given by
$$k_g = -\lambda \tau_g \tau_{g_1}$$

iii) if $x$ is a geodesic curve, then the geodesic curvature of the curve $x_1$ is given by
$$k_{g_1} = \lambda \tau_g \tau_{g_1}$$

By considering Liouville's formula, from Theorem 2 we can give the following corollary which gives the relations of the geodesic curvatures of parametric curves for the surfaces $S(u,v)$ and $S_1(u_1, v_1)$.

**Corollary 2.** *Let the pair $\{S(u,v), S_1(u_1,v_1)\}$ be a Bertrand pair surfaces with $\{x, x_1\}$ Bertrand $D$-pair. If the geodesic curvatures of parametric curves of $S(u,v)$ are $(k_g)_1$, $(k_g)_2$ and the geodesic curvatures of parametric curves of $S_1(u_1, v_1)$ are $(k_{g_1})_1$, $(k_{g_1})_2$ then we have*

$$\left(\frac{d\sigma}{ds} + (k_g)_1 \cos\sigma + (k_g)_2 \sin\sigma\right) - \left(\frac{d\sigma_1}{ds_1} + (k_{g_1})_1 \cos\sigma_1 + (k_{g_1})_2 \sin\sigma_1\right)$$
$$= \lambda\left[\left(\frac{d\sigma}{ds} + (k_g)_1 \cos\sigma + (k_g)_2 \sin\sigma\right)\left(\frac{d\sigma_1}{ds_1} + (k_{g_1})_1 \cos\sigma_1 + (k_{g_1})_2 \sin\sigma_1\right) - \tau_g \tau_{g_1}\right]$$

*where $\sigma$ is the angle between the curve $x(s)$ and the parametric curve $v = $ constant of $S(u,v)$ and $\sigma_1$ is the angle between the curve $x_1(s_1)$ and the parametric curve $v_1 = $ constant of $S_1(u_1, v_1)$.*

**Theorem 3.** *Let $\{x, x_1\}$ be Bertrand $D$-pair. Then the following relations hold:*

i) $k_{n_1} = k_n \dfrac{ds}{ds_1} + \dfrac{d\theta}{ds_1}$

ii) $\tau_g \dfrac{ds}{ds_1} = k_{g_1} \sin\theta + \tau_{g_1} \cos\theta$

iii) $k_g \dfrac{ds}{ds_1} = k_{g_1} \cos\theta + \tau_{g_1} \sin\theta$

iv) $\tau_{g_1} = (k_g \sin\theta + \tau_g \cos\theta)\dfrac{ds}{ds_1}$

**Proof: i)** By differentiating the equation $\langle T, T_1 \rangle = \cos\theta$ with respect to $s_1$ we have

$$\left\langle (k_g g + k_n n)\frac{ds}{ds_1}, T_1 \right\rangle + \left\langle T, k_{g_1} g_1 + k_{n_1} n_1 \right\rangle = -\sin\theta \frac{d\theta}{ds_1}.$$

Using the fact that the direction of $g_1$ coincides with the direction of $g$ and

$$T_1 = \cos\theta T + \sin\theta n,$$
$$n_1 = -\sin\theta T + \cos\theta n, \tag{27}$$

we easily get that

$$k_{n_1} = k_n \frac{ds}{ds_1} + \frac{d\theta}{ds_1}.$$

**ii)** By differentiating the equation $\langle n, g_1 \rangle = 0$ with respect to $s_1$ we have

$$\left\langle (-k_n T - \tau_g g)\frac{ds}{ds_1}, g_1 \right\rangle + \left\langle n, k_{g_1} T_1 + \tau_{g_1} n_1 \right\rangle = 0.$$

By (27) we obtain

$$\tau_g \frac{ds}{ds_1} = k_{g_1} \sin\theta + \tau_{g_1} \cos\theta$$

**iii)** By differentiating the equation $\langle T, g_1 \rangle = 0$ with respect to $s_1$ we get

$$\left\langle (k_g g + k_n n)\frac{ds}{ds_1}, g_1 \right\rangle + \left\langle T, (-k_{g_1} T_1 + \tau_{g_1} n_1) \right\rangle = 0.$$

From (27) it follows that

$$k_g \frac{ds}{ds_1} = k_{g_1} \cos\theta + \tau_{g_1} \sin\theta.$$

**iv)** By differentiating the equation $\langle n_1, g \rangle = 0$ with respect to $s_1$ we obtain

$$\left\langle -k_{n_1} T_1 - \tau_{g_1} g_1, g \right\rangle + \left\langle n_1, (-k_g T + \tau_g n)\frac{ds}{ds_1} \right\rangle = 0,$$

and using the fact that direction of $g_1$ coincides with the direction of $g$ and

$$T = \cos\theta T_1 - \sin\theta n_1,$$
$$n = \sin\theta T_1 + \cos\theta n_1,$$

we get

$$\tau_{g_1} = (k_g \sin\theta + \tau_g \cos\theta)\frac{ds}{ds_1}.$$

Let now $x$ be a Bertrand $D$-curve and $x_1$ be a Bertrand partner $D$-curve of $x$. From the first equation of (3) and by using the fact that $n_1 = -\sin\theta T + \cos\theta n$ we have

$$k_{g_1} = \left[(1 + \lambda k_g)\cos\theta - \lambda\tau_g \sin\theta\right]\left[k_g + \lambda k_g^2 + \lambda\tau_g^2\right]\left(\frac{ds}{ds_1}\right)^3. \tag{28}$$

Then we can give the following corollary.

***Corollary 3.*** *Let $x$ be a Bertrand $D$-curve and $x_1$ be a Bertrand partner $D$-curve of $x$. Then the relations between the geodesic curvature $k_{g_1}$ of $x_1(s_1)$ and the geodesic curvature $k_g$ and the geodesic torsion $\tau_g$ of $x(s)$ are given as follows.*

*i) If $x$ is a geodesic curve, then the geodesic curvature $k_{g_1}$ of $x_1(s_1)$ is*

$$k_{g_1} = \lambda \tau_g^2 (\cos\theta - \lambda\tau_g \sin\theta)\left(\frac{ds}{ds_1}\right)^3. \tag{29}$$

***ii)*** *If $x$ is a principal line, then the relation between the geodesic curvatures $k_{g_1}$ and $k_g$ is given by*

$$k_{g_1} = k_g(1+\lambda k_g)^2 \cos\theta \left(\frac{ds}{ds_1}\right)^3. \tag{30}$$

Similarly, From the second equation of (3) and by using the fact that $g$ is coincident with $g_1$, i.e., $n_1 = -\sin\theta T + \cos\theta n$, the geodesic torsion $\tau_{g_1}$ of $x_1$ is given by

$$\tau_{g_1} = \left[(\tau_g + \lambda k_g \tau_g)\cos^2\theta + (k_g + \lambda k_g^2 - \lambda\tau_g^2)\sin\theta\cos\theta - \lambda\tau_g k_g \sin^2\theta\right]\left(\frac{ds}{ds_1}\right)^2. \tag{31}$$

From (31) we can give the following corollary.

**Corollary 4.** *Let $x$ be a Bertrand $D$-curve and $x_1$ be a Bertrand partner $D$-curve of $x$. Then the relations between the geodesic torsion $\tau_{g_1}$ of $x_1(s_1)$ and the geodesic curvature $k_g$ and the geodesic torsion $\tau_g$ of $x(s)$ are given as follows.*

***i)*** *If $x$ is a geodesic curve then the geodesic torsion of $x_1$ is*

$$\tau_{g_1} = \tau_g \cos\theta\left[\cos\theta - \lambda\tau_g \sin\theta\right]\left(\frac{ds}{ds_1}\right)^2. \tag{32}$$

***ii)*** *If $x$ is a principal line then the relation between $\tau_{g_1}$ and $k_g$ is*

$$\tau_{g_1} = k_g(1+\lambda k_g)\sin\theta\cos\theta\left(\frac{ds}{ds_1}\right)^2. \tag{33}$$

Furthermore, by using (12) and (13), from (32) and (33) we have the following corollary.

**Corollary 5.** ***i)*** *Let $\{x, x_1\}$ be Bertrand $D$-pair and let $x$ be a geodesic line. Then the geodesic torsion $\tau_{g_1}$ of $x_1(s_1)$ is given by*

$$\tau_{g_1} = \tau_g(1-\lambda k_{g_1})\left[(1-\lambda k_{g_1}) + \lambda^2 \tau_g \tau_{g_1}\right]. \tag{34}$$

***ii)*** *Let $\{x, x_1\}$ be Bertrand $D$-pair and let $x$ be a principal line. Then the relation between the geodesic curvatures $k_g$ and $k_{g_1}$ is given as follows*

$$k_g(1+\lambda k_g)(1-\lambda k_{g_1}) = -\frac{1}{\lambda} = constant. \tag{35}$$

### 4. Conclusions
In this paper, the definition and characterizations of Bertrand partner $D$-curves are given which is a new study of associated curves lying on surfaces. It is shown that the definition and the characterizations of Bertrand partner $D$-curves include those of Bertrand partner curves in some special cases. Furthermore, the relations between the geodesic curvatures, the normal curvatures and the geodesic torsions of these curves are given.